\documentclass[12pt]{article}
\usepackage[utf8]{inputenc}
\usepackage[english]{babel}
\usepackage{amsmath,amssymb,amsthm,mathtools}
\usepackage{indentfirst}
\usepackage{tikz,float}
\usepackage[colorlinks=true,allcolors=blue]{hyperref}

\newtheorem{theorem}{Theorem}
\newtheorem{lemma}{Lemma}
\newtheorem{corollary}{Corollary}
\newcommand{\R}{\mathbb{R}}
\renewcommand{\P}{\mathbb{P}}
\newcommand{\eps}{\varepsilon}
\renewcommand{\le}{\leqslant}

\renewcommand{\leq}{\leqslant}
\renewcommand{\geq}{\geqslant}
\begin{document}
\title{Upper bounds for covering sets in $\ell_p^n$ by $\ell_p$-balls of prescribed diameter}
\author{A. I. Bikeev, V. A. Voronov, A. M. Raigorodskii}
\date{}
\maketitle

\begin{abstract}
    We obtain a new upper bound on the number of $\ell_p$-balls of diameter strictly less than $b$ needed to cover an arbitrary set of diameter $1$ in $\ell_p^n$, where $1\leq p<\infty$ and $0<b\leq1$.
\end{abstract}

\section{Introduction}
Borsuk's problem asks for the smallest number $f(n)$ such that every bounded set of positive diameter in Euclidean space $\R^n$ can be partitioned into at most $f(n)$ sets of strictly smaller diameter. Throughout the paper, asymptotic estimates are understood as $n\to\infty$, with all other parameters fixed. The asymptotic bounds are

\begin{equation}\label{eq:fn}
	\left(\left(\frac{2}{\sqrt{3}}\right)^{\sqrt{2}}+o(1)\right)^{\sqrt{n}} = (1.2255\ldots + o(1))^{\sqrt{n}} \le f(n) \le \left(\sqrt{\frac{3}{2}}+o(1)\right)^n = (1.224\ldots+o(1))^n.
\end{equation}

The lower bound was obtained by Raigorodskii \cite{Rai1}, and the upper bound was independently proved in \cite{Sch} and \cite{BL} in two different ways. Schramm \cite{Sch} proved this estimate using a covering of a constant-width body by its smaller positive homothetic copies. Bourgain and Lindenstrauss~\cite{BL} obtained the same exponential base by covering a set of diameter $1$ with balls of diameter less than $1$. Recently Akhiiarov, Koshelev, Kupavskii, and Raigorodskii~\cite{AKKR26} improved  the lower bound to $(1.3700\ldots+o(1))^{\sqrt{n}}$ and obtained $\limsup_{n\to\infty}f(n)^{1/\sqrt{n}}\geq 1.3839\ldots$

For $b\in(0,1]$, let $\chi(n,b)$ be the smallest integer such that every set of diameter $1$ in Euclidean space $\R^n$ can be covered by at most $\chi(n,b)$ sets of diameter strictly less than $b$. Define $c(n,b)$ in the same way, with the covering sets required to be balls. Replacing the strict diameter condition by a non-strict one in the definition of $c(n,b)$ does not change the exponential base: every Euclidean ball can be covered by $n+1$ balls of strictly smaller diameter, and $(n+1)^{1/n}\to1$. The result of Bourgain and Lindenstrauss immediately yields the estimate
\begin{equation}\label{eq:BL_standard}
\chi(n,1) \leq c(n,1) \leq \left(\sqrt{\frac{3}{2}}+o(1)\right)^n.
\end{equation}

In a recent paper \cite{BR26}, the authors provided a detailed survey of the upper bounds on the quantity $\chi(n,b)$ and proved the estimate
\begin{equation}\label{eq:BL_b}
    \chi(n,b) \leq c(n,b) \leq \left(\sqrt{\frac{1}{b^2} + \frac{1}{2}} + o(1) \right)^n.
\end{equation}
The bound (\ref{eq:BL_b}) is the strongest currently known upper bound on the quantities $\chi(n,b)$ and $c(n,b)$ for any $b \in (0,1)$. Clearly, for $b =1$, this bound coincides with (\ref{eq:BL_standard}).
See also \cite{BR2} and \cite{BVR} for some related results.

Let $\ell_p^n$ denote $\R^n$ equipped with the norm
 \[\|x\|_p = \left(|x_1|^p + |x_2| ^p + \dots + |x_n|^p \right)^{\frac{1}{p}}, \quad p\in[1,+\infty);\]
 \[\|x\|_{\infty} = \max\limits_{1 \leq i \leq n} |x_i|.\]

Define $\chi_p(n,b)$ and $c_p(n,b)$ analogously, using the $\ell_p$-metric and $\ell_p$-balls. Thus $\chi_2(n,b)=\chi(n,b)$ and $c_2(n,b)=c(n,b)$.

Let $N(K,L)$ denote the smallest number of translates of $L$ needed to cover $K$. Write $B^n_{p,R}$ for the closed ball of radius $R$ centered at the origin in $\ell_p^n$, and put $B^n_R=B^n_{2,R}$.

A classical result of Rogers (see \cite{Ro63, Ver, BW}) gives the asymptotic number of smaller Euclidean balls needed to cover a Euclidean ball. For $0<r<R$,
\begin{equation}\label{eq:rogers}
    N(B^n_R,B^n_r)=\left(\frac{R}{r}+o(1)\right)^n.
\end{equation}

Let $J_p^n$ be the smallest radius such that every set of diameter $1$ in $\ell_p^n$ is contained in an $\ell_p$-ball of that radius, and let $J_p=\sup_{n\in\mathbb N}J_p^n$. These are the Jung constants, named after Jung, who proved \cite{Jung} that
\begin{equation}\label{eq:Jung_classic}
    J_2^n=\sqrt{\frac{n}{2(n+1)}},\qquad J_2=\frac{1}{\sqrt2}.
\end{equation}
Pichugov \cite{Pichugov} obtained the following upper bounds for $1\leq p<\infty$:
\begin{equation}\label{eq:Pichugov}
    J_p\leq\begin{cases}
        2^{1/p-1}, & 1\leq p\leq2,\\
        2^{-1/p}, & 2\leq p<\infty.
    \end{cases}
\end{equation}

Chen, Gao, Li, and Wu \cite{CGLW} proved the following estimate, which applies to a discrete family of covering radii. For $1\leq p<\infty$ and every positive integer $k$,
\begin{equation}\label{eq:CGLW}
    N\left(B^n_{p,1},\left(\frac{n}{n+k}\right)^{1/p}B^n_{p,1}\right)
    \leq1+\sum_{i=1}^k2^i\binom{n}{i}\binom{k-1}{i-1}.
\end{equation}
Here $\binom{n}{i}=0$ when $i>n$.

For a Euclidean ball of diameter $1$, taking $p=2$ and $k=1$ in \eqref{eq:CGLW} gives a cover by at most $2n+1$ balls of diameter strictly less than $1$, which is stronger than the general bound \eqref{eq:BL_b} at $b=1$ for large $n$. For an arbitrary Euclidean set of diameter $1$, one can apply \eqref{eq:CGLW} to a containing Jung ball of radius $J_2^n$. By volume comparison, covering the entire Jung ball by balls of diameter less than $b$ requires at least $(2J_2^n/b)^n$ balls. The limiting base $\sqrt{2}/b$ exceeds $\sqrt{1/b^2+1/2}$ for every $b\in(0,1]$, so this approach gives a weaker bound than \eqref{eq:BL_b}.

The results of Rogers and Zong \cite{Rogers_density, RZ} imply that, for every centrally symmetric convex body $K\subset\R^n$ and every $b\in(0,1]$,
\begin{equation}\label{eq:RZ_classic}
    N(K,bK)\leq\left(1+\frac1b+o(1)\right)^n.
\end{equation}
Their method also yields
\begin{equation}\label{eq:rz_chi}
    \chi(n,b)=\chi_2(n,b)\leq\left(1+\frac1b+o(1)\right)^n,
    \qquad b\in(0,1].
\end{equation}
In the Euclidean case, \eqref{eq:BL_b} gives a strictly smaller exponential base for every $b\in(0,1]$, since $(1+1/b)^2-(1/b^2+1/2)=2/b+1/2>0$.

Combining a Jung ball with \eqref{eq:RZ_classic} gives a bound for coverings by $\ell_p$-balls. To respect the strict diameter condition in the definition of $c_p(n,b)$, first choose $0<b'<b$. Then
\[
    \chi_p(n,b)\leq c_p(n,b)
    \leq N\left(B^n_{p,J_p},B^n_{p,b'/2}\right)
    \leq\left(1+\frac{2J_p}{b'}+o(1)\right)^n.
\]
Taking the limsup of the $n$-th roots as $n\to\infty$, then letting $b'\uparrow b$ and using \eqref{eq:Pichugov}, we obtain
\begin{equation}\label{eq:RZ}
    \chi_p(n,b)\leq c_p(n,b)\leq
    \begin{cases}
        \left(\dfrac{2^{1/p}}b+1+o(1)\right)^n, & 1\leq p\leq2,\\[1ex]
        \left(\dfrac{2^{1-1/p}}b+1+o(1)\right)^n, & 2\leq p<\infty.
    \end{cases}
\end{equation}

For $p=\infty$, the exact values follow from the geometry of cubes. Every set of diameter $1$ is contained in a translate of $[0,1]^n$, which can be covered by a coordinate grid of smaller cubes. A grid of points with mutual $\ell_\infty$-distances at least $b$ gives the matching lower bound. Thus, for $0<b\leq1$,
\begin{equation}\label{eq:infty1}
    \chi_\infty(n,b)=c_\infty(n,b)
    =\inf_{0<\eps<b}N(B^n_{\infty,1},B^n_{\infty,b-\eps})
    =\left(\left\lfloor\frac1b\right\rfloor+1\right)^n.
\end{equation}
In particular,
\begin{equation}\label{eq:infty2}
    \chi_\infty(n,1)=c_\infty(n,1)=2^n.
\end{equation}

In this note, we obtain new upper bounds on $\chi_p(n,b)$ and $c_p(n,b)$. The standard approach leading to~\eqref{eq:RZ} places the set in a Jung ball and then covers the entire ball using the Rogers--Zong estimate. We use the diameter of the set itself to improve this bound. First, we choose random covering balls with centers in a Minkowski enlargement of its convex hull. The Brunn--Minkowski inequality controls the volume of this enlargement in terms of the diameter and yields the exponential base $1+1/b$ in every finite-dimensional normed space, independently of its Jung constant.

For $p$ close to $2$, a further improvement comes from transferring the Euclidean estimate \eqref{eq:BL_b} to $\ell_p^n$. After translating the set into a Jung ball, we apply the coordinatewise map $x_i\mapsto\operatorname{sgn}(x_i)|x_i|^{p/2}$, cover its Euclidean image, and transfer the covering back using the inverse map. The resulting bound depends on the diameter of the image and the distortion under the inverse map. We sharpen these estimates by integrating derivatives along chords of the containing ball, taking account of the smaller norms in the interior of each chord. For $p<2$, this improves the inverse-map estimate; for $p>2$, it improves the bound on the diameter of the Euclidean image. When $p$ is sufficiently close to $2$, the gain from the Euclidean covering estimate outweighs the distortion introduced by the maps. At $p=2$, the maps are the identity, and we recover \eqref{eq:BL_b}.

\section{Main results}

Our first theorem gives a covering bound in an arbitrary finite-dimensional normed space. Its proof uses a standard probabilistic argument, although we have not found this precise statement in the literature. Borsuk-type problems in Minkowski spaces have an extensive body of work of their own, see \cite{RS}.
\begin{theorem}\label{t:1b}
    Let $B$ be the unit ball of a norm on $\R^n$, where $n\geq2$, and let $b>0$.
    Every set $A$ of diameter at most $1$ in this norm can be covered by at most
    \begin{equation}\label{eq:t1}
        e^{2/(b(b+1))}\left(1+\frac1b\right)^n
        \bigl(n\log(1+4n)+1\bigr)
        =\left(1+\frac1b+o(1)\right)^n
    \end{equation}
    translates of $\frac b2B$.
\end{theorem}

The main theorem combines this bound with estimates obtained by transferring Euclidean coverings to $\ell_p^n$.
\begin{theorem}\label{t:closed}
    Let $1\leq p<\infty$ and $0<b\leq1$. Then
    \begin{equation}\label{eq:tmain_c}
        c_p(n,b)\leq\left(\min\left\{1+\frac1b,x(p,b)\right\}+o(1)\right)^n,
    \end{equation}
    \begin{equation}\label{eq:tmain_chi}
        \chi_p(n,b)\leq\left(\min\left\{1+\frac1b,y(p,b)\right\}+o(1)\right)^n,
    \end{equation}
    where
    \begin{equation}\label{eq:x}
        x(p,b)=\begin{cases}
\sqrt{\dfrac12+\dfrac{2^{1+2/p}}{p^2b^2}\left(1-\dfrac{p^2b^2}{24\cdot2^{2/p}}\right)^{(2-p)/p}}, & 1\leq p\leq2,\\[3ex]
\sqrt{\dfrac12+\dfrac{p^2 2^{2p-7+2/p}}{b^p}\left(1-\dfrac1{2^p-1}\right)^{(p-2)/p}}, & 2\leq p<\infty,
        \end{cases}
    \end{equation}
    and
    \begin{equation}\label{eq:y}
        y(p,b)=\begin{cases}
\sqrt{\dfrac12+\dfrac{2^{1+2/p}}{p^2b^2}\left(1-\dfrac{p^2b^2}{6\cdot2^{2/p}}\right)^{(2-p)/p}}, & 1\leq p\leq2,\\[3ex]
\sqrt{\dfrac12+\dfrac{p^2 2^{p-5+2/p}}{b^p}\left(1-\dfrac1{2^p-1}\right)^{(p-2)/p}}, & 2\leq p<\infty.
        \end{cases}
    \end{equation}
\end{theorem}

\begin{figure}[htbp]
    \centering
    \begin{tabular}{cc}
        \includegraphics[width=0.45\textwidth]{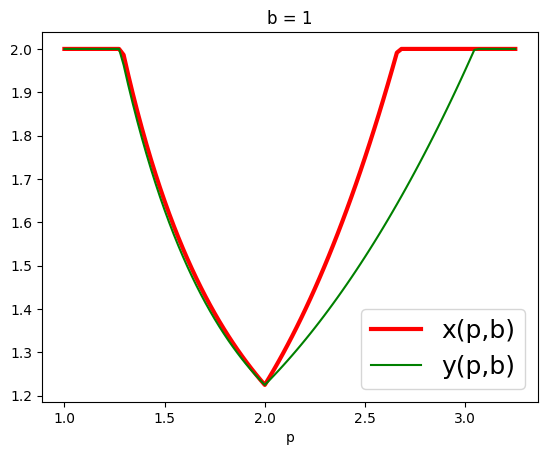} &
        \includegraphics[width=0.45\textwidth]{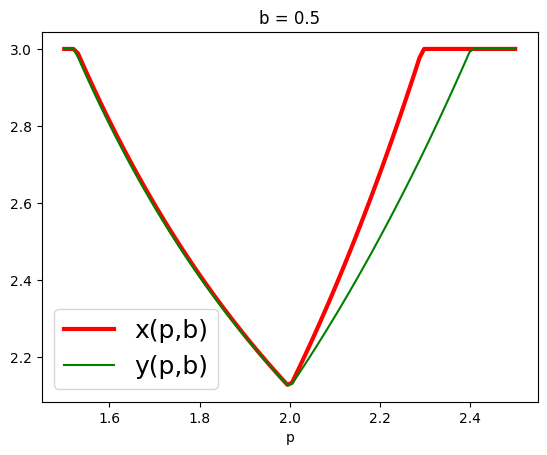}
        \\[2mm]
        \includegraphics[width=0.45\textwidth]{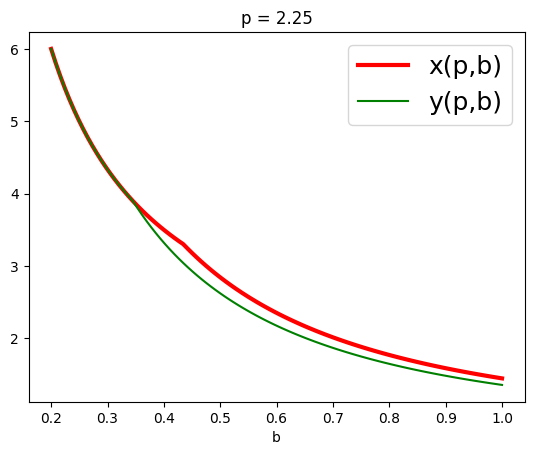}
        &
        \includegraphics[width=0.45\textwidth]{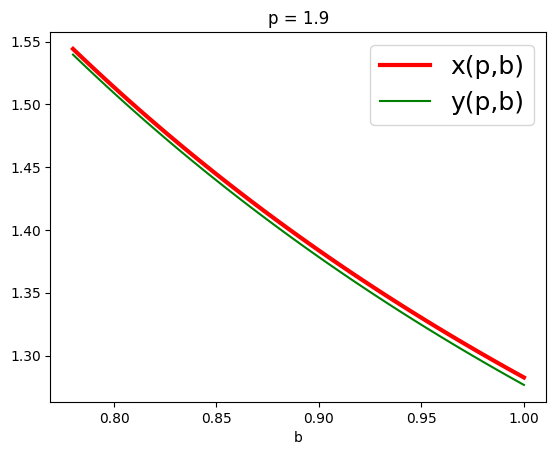}
    \end{tabular}
    \caption{Exponential bases in the upper bounds in Theorem \ref{t:closed}}
    \label{fig:xy}
\end{figure}

The bounds in Theorem~\ref{t:closed} improve on \eqref{eq:RZ}. For each fixed $b$, both $x(p,b)$ and $y(p,b)$ are smaller than $1+1/b$ when $p$ is sufficiently close to $2$. The formulas give $y(p,b)\leq x(p,b)$ throughout the stated range. For $p\geq2$, the difference is quantified by
\[
    x(p,b)^2-\frac12=2^{p-2}\left(y(p,b)^2-\frac12\right).
\]
The smaller base for $\chi_p(n,b)$ comes from allowing the covering sets to be inverse images of Euclidean balls, with their diameters controlled directly. For $c_p(n,b)$, the inverse images must fit inside $\ell_p$-balls of the required radius. Figure~\ref{fig:xy} compares the resulting bases. At $p=2$, both expressions reduce to $\sqrt{1/2+1/b^2}$, the base in \eqref{eq:BL_b}.

\section{Proof of Theorem \ref{t:1b}}
\begin{proof}
We may assume that $A$ is nonempty and replace it by its closed convex hull $K$, which is compact and still has diameter at most $1$. The argument first covers a finite net of $K$ by randomly chosen balls and then enlarges the balls to cover $K$.

Choose the net spacing so that the initial covering radius remains positive for every $b>0$:
\[
    a=\min\{b,1\},\qquad \eta=\frac{a}{2n},\qquad
    q=\frac{b-a/n}{2}>0,\qquad T=K+qB.
\]
Write $X-Y=\{x-y:x\in X,\ y\in Y\}$ and let $\operatorname{vol}$ denote $n$-dimensional volume. Since $T-T\subset(1+2q)B$, the Brunn--Minkowski inequality (see \cite{Melnikov}) gives
\[
    2^n\operatorname{vol}T
    \leq\operatorname{vol}(T-T)
    \leq(1+2q)^n\operatorname{vol}B.
\]

Let $M\subset K$ be a maximal set whose distinct points have distance greater than $\eta$. Then $M$ is an $\eta$-net of $K$. The balls $x+\frac\eta2B$, $x\in M$, have disjoint interiors and lie in a ball of radius $1+\eta/2$ centered at any fixed point of $K$. Comparing volumes, we obtain
\[
    |M|\leq\left(1+\frac2\eta\right)^n
    =\left(1+\frac{4n}{a}\right)^n.
\]

The volume estimate for $T$ gives a uniform lower bound on the probability that a random ball covers a fixed point. Indeed, if $z$ is uniformly distributed in $T$, then $x+qB\subset T$ for every $x\in K$, and symmetry of $B$ yields
\[
    \P\{x\in z+qB\}
    =\frac{q^n\operatorname{vol}B}{\operatorname{vol}T}
    \geq\left(\frac{2q}{1+2q}\right)^n=:\pi_n.
\]
Choose $N=\lfloor\pi_n^{-1}\log|M|\rfloor+1$ independent uniform centers in $T$. The probability that some point of $M$ is uncovered is at most
\[
    |M|(1-\pi_n)^N\leq |M|e^{-N\pi_n}<1.
\]
Thus some choice of $N$ balls of radius $q$ covers $M$. Enlarging their radii to $q+\eta=b/2$ covers all of $K$.

It remains to bound $N$. Since $a\leq b$ and $n\geq2$,
\[
\begin{aligned}
    \pi_n^{-1}
    &=\left(1+\frac1{b-a/n}\right)^n\\
    &=\left(1+\frac1b\right)^n
      \left(1+\frac{a}{(nb-a)(b+1)}\right)^n
    \leq e^{2a/(b(b+1))}\left(1+\frac1b\right)^n.
\end{aligned}
\]
Consequently,
\[
    N\leq\pi_n^{-1}(\log|M|+1)
    \leq e^{2a/(b(b+1))}\left(1+\frac1b\right)^n
       \left(n\log\left(1+\frac{4n}{a}\right)+1\right).
\]
For $a=1$, this is the required bound. For $a=b<1$, use
\[
    \log(1+4n/b)\leq\frac1b\log(1+4n),\qquad
    \log\frac1b\leq\frac{1-b}{b}\leq\frac{2(1-b)}{b(b+1)}.
\]
The first inequality bounds the parenthetical factor by $b^{-1}(n\log(1+4n)+1)$, and the second absorbs $b^{-1}$ into the exponential prefactor. This proves \eqref{eq:t1}.
\end{proof}

\section{Notation and auxiliary results}

We transfer coverings between normed spaces using coordinatewise power maps. For $q_1,q_2\in[1,\infty)$, define $\zeta_{q_1,q_2}:\ell_{q_1}^n\to\ell_{q_2}^n$ (known as the Mazur map) by
\[
    \zeta_{q_1,q_2}(x)_i=\operatorname{sgn}(x_i)|x_i|^{q_1/q_2},
    \qquad 1\leq i\leq n,
\]
where $\operatorname{sgn}(0)=0$. The next lemma records the continuity estimates needed to control the diameters of images and inverse images.

\begin{lemma}\label{l:main}
Let $p_1,p_2\in[1,\infty)$ and $\alpha=p_1/p_2$.
\begin{enumerate}
    \item The map $\zeta_{p_1,p_2}$ is a homeomorphism, with inverse $\zeta_{p_2,p_1}$.
    \item For every $x\in\ell_{p_1}^n$,
    \[
        \|\zeta_{p_1,p_2}(x)\|_{p_2}=\|x\|_{p_1}^\alpha.
    \]
    In particular, $\zeta_{p_1,p_2}(B^n_{p_1,r})=B^n_{p_2,r^\alpha}$.
    \item If $\alpha\leq1$, then, for all $x,y\in\ell_{p_1}^n$,
    \[
        \|\zeta_{p_1,p_2}(x)-\zeta_{p_1,p_2}(y)\|_{p_2}
        \leq2^{1-\alpha}\|x-y\|_{p_1}^\alpha.
    \]
    \item If $\alpha\geq1$, $R>0$, and $\|x\|_{p_1},\|y\|_{p_1}\leq R$, then
    \[
        \|\zeta_{p_1,p_2}(x)-\zeta_{p_1,p_2}(y)\|_{p_2}
        \leq\alpha R^{\alpha-1}\|x-y\|_{p_1}.
    \]
    For $\alpha>1$, the map is continuously differentiable, and its Jacobian at $z$, viewed as an operator from $\ell_{p_1}^n$ to $\ell_{p_2}^n$, has norm at most $\alpha\|z\|_{p_1}^{\alpha-1}$. For $\alpha=1$, the map is the identity.
\end{enumerate}
\end{lemma}

The third and fourth statements are the standard H\"older and local-Lipschitz estimates for the Mazur map; see \cite{BeLi}.

\begin{proof}
The inverse formula and the norm identity follow directly from the definition.

For the third assertion, first consider one coordinate. If $0<\alpha\leq1$ and $s,t\geq0$, concavity gives
\[
    |s^\alpha-t^\alpha|\leq|s-t|^\alpha,
    \qquad s^\alpha+t^\alpha\leq2^{1-\alpha}(s+t)^\alpha.
\]
Use the first inequality when the coordinates have the same sign and the second when they have opposite signs. Raising to the power $p_2$ gives
\[
    |\zeta_{p_1,p_2}(x)_i-\zeta_{p_1,p_2}(y)_i|^{p_2}
    \leq\begin{cases}
        |x_i-y_i|^{p_1}, & x_i y_i\geq0,\\
        2^{p_2-p_1}|x_i-y_i|^{p_1}, & x_i y_i<0.
    \end{cases}
\]
Summing over the coordinates and taking the $p_2$-th root proves the assertion.

For the fourth assertion, the case $\alpha=1$ is immediate. If $\alpha>1$, the Jacobian $J(z)$ is diagonal, with
\[
    (J(z)h)_i=\alpha|z_i|^{\alpha-1}h_i.
\]
H\"older's inequality, with conjugate exponents $\alpha/(\alpha-1)$ and $\alpha$, yields
\[
    \sum_{i=1}^n|z_i|^{(\alpha-1)p_2}|h_i|^{p_2}
    \leq\left(\sum_{i=1}^n|z_i|^{p_1}\right)^{(\alpha-1)/\alpha}
         \left(\sum_{i=1}^n|h_i|^{p_1}\right)^{1/\alpha}.
\]
Taking the $p_2$-th root proves the operator norm bound. Integrating this bound along the segment from $x$ to $y$, which lies in $B^n_{p_1,R}$, gives the required distance estimate.
\end{proof}

To improve the derivative estimate along a segment, we quantify how far its interior lies inside a norm ball. For $p\in[1,\infty)$ and $t\in[0,1]$, define
\[
    \xi_p(t)=\inf_{s\in\R}
    \left\{(1-t)|s-t|^p+t|s+1-t|^p-|s|^p\right\}.
\]

\begin{lemma}
For every $p\in[1,\infty)$, $t\in[0,1]$, and $x,y\in\ell_p^n$, we have $\xi_p(t)\geq0$ and
\begin{equation}\label{eq:txy}
    \|(1-t)x+ty\|_p^p
    \leq(1-t)\|x\|_p^p+t\|y\|_p^p-\xi_p(t)\|x-y\|_p^p.
\end{equation}
Moreover,
\begin{equation}\label{eq:xi2}
    \xi_2(t)=t(1-t),
\end{equation}
\begin{equation}\label{eq:xipt}
    \xi_p(t)\leq2^{-p}.
\end{equation}
\end{lemma}

\begin{proof}
Convexity of $s\mapsto|s|^p$ gives $\xi_p(t)\geq0$. To obtain \eqref{eq:txy}, apply the definition to each coordinate: when $x_i\ne y_i$, set $s=((1-t)x_i+ty_i)/(y_i-x_i)$ and multiply by $|y_i-x_i|^p$. Coordinates with $x_i=y_i$ contribute equality. Summing proves the claim.

For $p=2$, expanding the squares shows that the expression defining $\xi_2(t)$ equals $t(1-t)$ independently of $s$. Finally, substituting $s=t-1/2$ gives
\[
    \xi_p(t)\leq2^{-p}-|t-1/2|^p\leq2^{-p},
\]
which proves \eqref{eq:xipt}. For each fixed $s$, the expression defining $\xi_p(t)$ is continuous in $t$. Its infimum over $s$ is therefore upper semicontinuous, so $\xi_p$ is measurable.
\end{proof}

\begin{corollary}\label{cor:txy}
If $\|x\|_p,\|y\|_p\leq r$ and $t\in[0,1]$, then
\[
    \|(1-t)x+ty\|_p^p\leq r^p-\xi_p(t)\|x-y\|_p^p.
\]
\end{corollary}

We now define the parameters in an intermediate covering bound. For $1\leq p\leq2$, put $D_p=2^{1-p/2}$ and, for $0\leq r\leq D_p$, let
\[
    I_p(r)=\frac{2r}{p}\int_0^1
       \left(2^{1-p}-t(1-t)r^2\right)^{1/p-1/2}\,dt.
\]
For $0<b\leq1$, define $q_p(b)\in(0,D_p/2]$ and $d_p(b)\in(0,D_p]$ by
\[
    I_p(q_p(b))=b/2,\qquad I_p(d_p(b))=b.
\]
Their existence and uniqueness are established below.

For $2\leq p<\infty$, define the following quantities, with the integral understood in the Lebesgue sense:
\[
    E_p=\frac p2\int_0^1\left(\frac12-\xi_p(t)\right)^{1/2-1/p}\,dt,
    \qquad\delta_p=\min\{\sqrt2,E_p\},
\]
\[
    q'_p(b)=\min\{\delta_p/2,2^{1-p}b^{p/2}\}.
\]
The integral is well defined because $\xi_p$ is measurable and $0\leq\xi_p(t)\leq2^{-p}\leq1/4$. The parameter $\delta_p$ bounds the diameter of the Euclidean image, while $q'_p(b)$ is a covering radius that gives $\ell_p$-balls of radius at most $b/2$ after applying the inverse map.

Define
\[
    \alpha(p,b)=\begin{cases}
        \sqrt{\dfrac12+\dfrac{2^{2-p}}{4q_p(b)^2}}, & 1\leq p\leq2,\\[3ex]
        \sqrt{\dfrac12+\dfrac{\delta_p^2}{4q'_p(b)^2}}, & 2\leq p<\infty,
    \end{cases}
\]
and
\[
    \beta(p,b)=\begin{cases}
        \sqrt{\dfrac12+\dfrac{2^{2-p}}{d_p(b)^2}}, & 1\leq p\leq2,\\[3ex]
        \sqrt{\dfrac12+2^{p-2}\delta_p^2b^{-p}}, & 2\leq p<\infty.
    \end{cases}
\]

\begin{lemma}\label{l:great}
For $1\leq p<\infty$ and $0<b\leq1$,
\begin{equation}\label{eq:lmain_c}
    c_p(n,b)\leq\left(\min\left\{1+\frac1b,\alpha(p,b)\right\}+o(1)\right)^n,
\end{equation}
\begin{equation}\label{eq:lmain_chi}
    \chi_p(n,b)\leq\left(\min\left\{1+\frac1b,\beta(p,b)\right\}+o(1)\right)^n.
\end{equation}
\end{lemma}

This lemma gives bounds at least as strong as those in Theorem~\ref{t:closed}, but its parameters involve integrals and inverse functions. Theorem~\ref{t:closed} follows by estimating these parameters explicitly.

\section{Proof of Lemma \ref{l:great}}

The base $1+1/b$ follows from Theorem~\ref{t:1b}, applied first with a parameter $b'<b$ and then with $b'\uparrow b$. We therefore prove the bounds with bases $\alpha(p,b)$ and $\beta(p,b)$.

We use two elementary observations about covering bounds. First, a nonempty set of diameter at most $D>0$ can be enlarged to one of diameter exactly $D$: after taking its compact closure, add a point $z$ for which $\max_x\|z-x\|_2=D$. Such a point exists by continuity, starting at a point of the set and moving to infinity. Thus \eqref{eq:BL_b}, after scaling, applies whenever $D$ is an upper bound on the image diameter.

Second, the constructions below initially give covering sets of diameter at most $b$. To obtain the strict condition, fix $b'<b$ and perform the construction with $b'$. If its exponential base is $F(p,b')$, then
\[
    \limsup_{n\to\infty}c_p(n,b)^{1/n}\leq F(p,b'),
\]
and the same reasoning applies to $\chi_p(n,b)$. Letting $b'\uparrow b$ then gives the base $F(p,b)$ by continuity. This order of limits requires no uniformity of the $o(1)$ term in $b'$.

\subsection{The case \texorpdfstring{$1\leq p<2$}{1<=p<2}}
Put
\[
    R=2^{1/p-1},\qquad S=R^{p/2}=2^{1/2-p/2},\qquad
    D=\sqrt2S=2^{1-p/2},
\]
\[
    f=\zeta_{p,2},\qquad g=\zeta_{2,p}=f^{-1}.
\]
Let $K\subset\ell_p^n$ have diameter $1$. By \eqref{eq:Pichugov}, we may translate $K$ into $B^n_{p,R}$. Lemma~\ref{l:main} then gives
\[
    L:=f(K)\subset B^n_{2,S},\qquad \operatorname{diam}_2L\leq D.
\]

We first estimate how much $g$ can increase a short Euclidean distance inside $B^n_{2,S}$. For $u,v\in B^n_{2,S}$, set $r=\|u-v\|_2$ and $z_t=(1-t)u+tv$. The Jacobian estimate in Lemma~\ref{l:main}, together with Corollary~\ref{cor:txy} for the Euclidean norm, gives
\[
    \|Jg(z_t)\|_{\ell_2\to\ell_p}
    \leq\frac2p\|z_t\|_2^{2/p-1}
    \leq\frac2p\left(S^2-t(1-t)r^2\right)^{1/p-1/2}.
\]
Integrating along the segment yields, whenever $r\leq D$,
\begin{equation}\label{eq:int}
    \|g(u)-g(v)\|_p
    \leq\frac{2r}{p}\int_0^1
       \left(S^2-t(1-t)r^2\right)^{1/p-1/2}\,dt
    =I_p(r).
\end{equation}

To choose covering radii using \eqref{eq:int}, we need $I_p$ to be invertible on the relevant interval. The derivative of its integrand, including the factor $r$, has the sign of
\[
    S^2-\frac2p t(1-t)r^2
    \geq S^2(1-1/p)\geq0,\qquad 0\leq r\leq D.
\]
For $r<D$, the derivative is positive, including when $p=1$. Hence $I_p$ is strictly increasing on $[0,D]$. Also, $I_p(0)=0$, and the inequalities $t(1-t)\leq1/4$ and $D^2=2S^2$ imply
\[
    I_p(D)\geq\frac{2D}{p}(S^2/2)^{1/p-1/2}
    =\frac2p>1,
\]
\[
    I_p(D/2)\geq\frac Dp(7S^2/8)^{1/p-1/2}
    =\frac1p\left(\frac74\right)^{1/p-1/2}>\frac12.
\]
Thus $q_p(b)$ and $d_p(b)$ exist uniquely and are continuous and strictly increasing in $b$.

We now transfer a Euclidean ball covering of $L$ to an $\ell_p$-ball covering of $K$. Applying \eqref{eq:BL_b} after scaling by $D$, we cover $L$ by $k$ Euclidean balls of radius $q_p(b)$, where
\[
    k\leq\left(\sqrt{\frac12+\frac{D^2}{4q_p(b)^2}}+o(1)\right)^n.
\]
We may project their centers $u_1,\ldots,u_k$ onto $B^n_{2,S}$: Euclidean projection does not increase the distance to any point of that ball. If $v\in L$ lies in the ball centered at $u_i$, then \eqref{eq:int} gives
\[
    \|g(u_i)-g(v)\|_p\leq I_p(q_p(b))=b/2.
\]
Hence the $\ell_p$-balls of radius $b/2$ centered at $g(u_i)$ cover $K$. This gives the base $\alpha(p,b)$.

For the bound on $\chi_p(n,b)$, apply \eqref{eq:BL_b} to cover $L$ by $m$ Euclidean balls $C_1,\ldots,C_m$ of diameter $d_p(b)$, with
\[
    m\leq\left(\sqrt{\frac12+\frac{D^2}{d_p(b)^2}}+o(1)\right)^n.
\]
Intersect each ball with $B^n_{2,S}$ before applying $g$, so that \eqref{eq:int} remains applicable. The sets $g(C_j\cap B^n_{2,S})$ cover $K$, and each has diameter at most $I_p(d_p(b))=b$. This gives the base $\beta(p,b)$.

\subsection{The case \texorpdfstring{$p>2$}{p>2}}
Put
\[
    R=2^{-1/p},\qquad S=R^{p/2}=\frac1{\sqrt2},\qquad
    f=\zeta_{p,2},\qquad g=\zeta_{2,p}=f^{-1}.
\]
As above, translate $K$ into $B^n_{p,R}$ and set $L=f(K)$. Then $L\subset B^n_{2,S}$, so $\operatorname{diam}_2L\leq\sqrt2$.

A second diameter bound follows by integrating the derivative of $f$ along segments in $B^n_{p,R}$. For $x,y\in K$ and $r=\|x-y\|_p\leq1$, Lemma~\ref{l:main} and Corollary~\ref{cor:txy} give
\[
    \|f(x)-f(y)\|_2
    \leq\frac{pr}{2}\int_0^1
       \left(\frac12-\xi_p(t)r^p\right)^{1/2-1/p}\,dt.
\]
For each fixed $t$, the integrand including the factor $r$ is increasing for $0\leq r\leq1$: its derivative has the sign of
\[
    \frac12-\frac p2\xi_p(t)r^p
    \geq\frac12-p2^{-p-1}>0.
\]
We may therefore set $r=1$ on the right-hand side and conclude that
\[
    \operatorname{diam}_2L\leq E_p,
    \qquad\operatorname{diam}_2L\leq\delta_p.
\]
The preceding estimate holds for every pair in $B^n_{p,R}$ at distance at most $1$. In dimension at least $2$, take $x=Re_1$ and $y=Re_2$. Then $\|x-y\|_p=2^{1/p}R=1$ and $\|f(x)-f(y)\|_2=\sqrt2S=1$, so $E_p\geq1$. Since $E_p$ is independent of dimension, this proves $\delta_p\geq1$.

The inverse map is controlled globally by the third assertion of Lemma~\ref{l:main}:
\[
    \|g(u)-g(v)\|_p\leq2^{1-2/p}\|u-v\|_2^{2/p}.
\]
Consequently, a Euclidean ball of radius $q'_p(b)$ centered at $u$ maps into an $\ell_p$-ball of radius at most $b/2$ centered at $g(u)$. Applying \eqref{eq:BL_b} to $L$ gives a covering with at most
\[
    \left(\sqrt{\frac12+\frac{\delta_p^2}{4q'_p(b)^2}}+o(1)\right)^n
\]
balls, which proves the bound with base $\alpha(p,b)$.

For arbitrary covering sets, the same inverse-map estimate allows Euclidean balls of diameter $2^{1-p/2}b^{p/2}$. This diameter is at most $1\leq\delta_p$, so \eqref{eq:BL_b} applies after scaling by $\delta_p$. Their images have $\ell_p$-diameter at most $b$, and the resulting exponential base is
\[
    \sqrt{\frac12+\frac{\delta_p^2}{2^{2-p}b^p}}=\beta(p,b).
\]

\subsection{The case \texorpdfstring{$p=2$}{p=2}}
Here $I_2(r)=r$, $q_2(b)=b/2$, $d_2(b)=b$, $\delta_2=1$, and $q'_2(b)=b/2$. Both formulas for each base agree and give
\[
    \alpha(2,b)=\beta(2,b)=\sqrt{\frac12+\frac1{b^2}}.
\]
Thus this case is exactly \eqref{eq:BL_b}, completing the proof of Lemma~\ref{l:great}.

\section{Proof of Theorem \ref{t:closed}}

\subsection{The case \texorpdfstring{$1\leq p\leq2$}{1<=p<=2}}
Use the same notation as in the proof of Lemma~\ref{l:great}:
\[
    R=2^{1/p-1},\qquad S=R^{p/2},\qquad D=\sqrt2S.
\]
We seek explicit lower bounds on the inverse radii $q_p(b)$ and $d_p(b)$. Since $x\mapsto x^{1/p-1/2}$ is concave and $\int_0^1t(1-t)\,dt=1/6$, Jensen's inequality gives, for $0\leq r\leq D$,
\begin{equation}\label{eq:Ipr}
    I_p(r)\leq\frac{2r}{p}\left(S^2-\frac{r^2}{6}\right)^{1/p-1/2}.
\end{equation}
The expression in parentheses is positive throughout this interval, and equality holds when $p=2$.

The simpler estimate obtained by replacing the parentheses with $S^2$ suggests an initial radius $d_0(b)$. We enlarge it to account for the negative quadratic term:
\[
    d_0(b)=\frac{pb}{2R^{1-p/2}},\qquad
    \widehat d(b)=d_0(b)\left(1-\frac{p^2b^2}{24R^2}\right)^{1/2-1/p}.
\]
First, we check that this enlarged radius lies in the interval where \eqref{eq:Ipr} applies. Since $R\geq1/\sqrt2$, $p\leq2$, and $b\leq1$,
\[
    0\leq\frac{p^2b^2}{24R^2}\leq\frac13,
    \qquad \frac{d_0(b)}D=\frac{pb}{2\sqrt2R}\leq\frac p2.
\]
Therefore,
\[
    \frac{\widehat d(b)}D\leq\frac p2\left(\frac32\right)^{1/p-1/2}\leq1.
\]
For the last inequality, use $\log x\leq x-1$ to obtain
\[
\begin{aligned}
    \log\frac p2+\left(\frac1p-\frac12\right)\log\frac32
    &\leq\frac{p-2}{2}+\frac{2-p}{2p}\log\frac32\\
    &=\frac{2-p}{2}\left(\frac{\log(3/2)}p-1\right)\leq0.
\end{aligned}
\]

We can now verify that $\widehat d(b)$ gives an admissible inverse radius. Since $d_0(b)\leq\widehat d(b)\leq D$, \eqref{eq:Ipr} implies
\[
\begin{aligned}
    I_p(\widehat d(b))
    &\leq\frac2p\widehat d(b)
       \left(S^2-\frac{\widehat d(b)^2}{6}\right)^{1/p-1/2}\\
    &\leq\frac2p\widehat d(b)
       \left(S^2-\frac{d_0(b)^2}{6}\right)^{1/p-1/2}\\
    &=\frac2p d_0(b)S^{2/p-1}=b.
\end{aligned}
\]
In the equality, the correction factor in $\widehat d(b)$ cancels because $d_0(b)^2/(6S^2)=p^2b^2/(24R^2)$. Strict monotonicity of $I_p$ now gives
\[
    q_p(b)\geq\widehat d(b/2),\qquad d_p(b)\geq\widehat d(b).
\]
Substituting these lower bounds into the definitions of the bases, we obtain
\[
    \alpha(p,b)\leq\sqrt{\frac12+\frac{2^{2-p}}{4\widehat d(b/2)^2}},
    \qquad
    \beta(p,b)\leq\sqrt{\frac12+\frac{2^{2-p}}{\widehat d(b)^2}}.
\]
Finally, substituting $R=2^{1/p-1}$ and simplifying shows that the right-hand sides are exactly $x(p,b)$ and $y(p,b)$ from \eqref{eq:x} and \eqref{eq:y}.

\subsection{The case \texorpdfstring{$p\geq2$}{p>=2}}
It remains to obtain an explicit upper bound on $\delta_p$. A lower bound on $\xi_p$ gives an upper bound on $E_p$, so we first estimate $\xi_p$ using repeated midpoint interpolation.

Define
\begin{equation}\label{eq:hpt}
    h_p(t)=\sum_{k=0}^\infty2^{1-(k+1)p}
       \operatorname{dist}(2^kt,\mathbb Z),\qquad t\in[0,1].
\end{equation}
The series converges uniformly, since its $k$-th term is at most $2^{-(k+1)p}$. Fix $s\in\R$ and $t\in[0,1]$, and put
\[
    \lambda(z)=|s-t+z|^p,\qquad z\in[0,1].
\]
The scalar form of Clarkson's inequality (see \cite{Clarkson}) gives, for $0\leq a\leq b\leq1$,
\[
    \lambda\left(\frac{a+b}{2}\right)
    \leq\frac{\lambda(a)+\lambda(b)}2-\left(\frac{b-a}{2}\right)^p.
\]
Let $L_k$ be the piecewise linear interpolant of $\lambda$ on the dyadic grid of spacing $2^{-k}$. On each grid interval, $L_k-L_{k+1}$ vanishes at the endpoints and has value at least $2^{-(k+1)p}$ at the midpoint. Its linearity on each half-interval therefore yields
\[
    L_k(z)-L_{k+1}(z)
    \geq2^{1-(k+1)p}\operatorname{dist}(2^kz,\mathbb Z).
\]
Summing over $k$ and using the uniform convergence $L_k\to\lambda$, we find
\[
    L_0(t)-\lambda(t)
    =(1-t)|s-t|^p+t|s+1-t|^p-|s|^p\geq h_p(t).
\]
Taking the infimum over $s$ proves
\begin{equation}\label{eq:eq}
    \xi_p(t)\geq h_p(t).
\end{equation}

The advantage of $h_p$ is that its integral is explicit. Each function $\operatorname{dist}(2^kt,\mathbb Z)$ has mean $1/4$ on $[0,1]$, so termwise integration gives
\begin{equation}\label{eq:inthp}
    \int_0^1h_p(t)\,dt
    =\frac14\sum_{k=0}^\infty2^{1-(k+1)p}
    =\frac1{2(2^p-1)}.
\end{equation}
By \eqref{eq:xipt}, $0\leq h_p(t)\leq\xi_p(t)\leq2^{-p}\leq1/4$. We may therefore apply \eqref{eq:eq} and then Jensen's inequality to the concave function $x\mapsto x^{1/2-1/p}$:
\[
\begin{aligned}
    E_p
    &\leq\frac p2\int_0^1\left(\frac12-h_p(t)\right)^{1/2-1/p}\,dt\\
    &\leq\frac p2\left(\frac12-\int_0^1h_p(t)\,dt\right)^{1/2-1/p}\\
    &=\frac p2\left(\frac12\left(1-\frac1{2^p-1}\right)\right)^{1/2-1/p}.
\end{aligned}
\]
Squaring gives the required explicit diameter bound:
\begin{equation}\label{eq:delta}
    \delta_p^2=\min\{2,E_p^2\}
    \leq p^2 2^{-3+2/p}
       \left(1-\frac1{2^p-1}\right)^{1-2/p}.
\end{equation}

As shown in the proof of Lemma~\ref{l:great}, $\delta_p\geq1$. Since $2^{2-p}b^{p/2}\leq1$, the minimum defining $q'_p(b)$ equals $2^{1-p}b^{p/2}$. Thus
\[
    \alpha(p,b)^2=\frac12+2^{2p-4}\delta_p^2b^{-p},
    \qquad
    \beta(p,b)^2=\frac12+2^{p-2}\delta_p^2b^{-p}.
\]
Substituting \eqref{eq:delta} into these identities gives $\alpha(p,b)\leq x(p,b)$ and $\beta(p,b)\leq y(p,b)$. Lemma~\ref{l:great} now proves Theorem~\ref{t:closed}.

\section{Statement of AI use}
We used GPT-6 Astra to obtain the main steps of the proofs of Theorems~\ref{t:1b} and~\ref{t:closed}. The authors independently verified the proofs and simplified some of the steps.

\clearpage
\section{Comparison of the bounds}

\begin{thebibliography}{99}

\bibitem{Rai1} A. M. Raigorodskii, \emph{On a bound in Borsuk's problem}, Russian Mathematical Surveys, 54 (1999), no. 2, 453--454.

\bibitem{Sch} O. Schramm, \emph{Illuminating sets of constant width}, Mathematika, 35 (1988), 180--189.

\bibitem{BL} J. Bourgain, J. Lindenstrauss, \emph{On covering a set in $\mathbb R^d$ by balls of the same diameter}, Geometric Aspects of Functional Analysis (J. Lindenstrauss and V. Milman, eds.), Lecture Notes in Math., 1469, Springer-Verlag, Berlin, 1991, 138--144.

\bibitem{AKKR26} A. Akhiiarov, M. Koshelev, A. Kupavskii, A. Raigorodskii, \emph{A polynomial method over formal power series}, preprint, \href{https://arxiv.org/abs/2610.39460}{arXiv:2610.39460}, 2026.

\bibitem{BR26} A. I. Bikeev, A. M. Raigorodskii, \href{https://doi.org/10.1134/S0001434626601899}{\emph{On upper bounds for the number of parts in the problem of partitioning sets into parts of smaller diameter}}, Math. Notes 119 (2026), 406--415.

\bibitem{BR2} A. I. Bikeev, A. M. Raigorodskii, \emph{On partitioning sets in $\mathbb{R}^n$ into parts that avoid a specified distance}, accepted for publication in Math. Notes.

\bibitem{BVR} A. I. Bikeev, V. A. Voronov, A. M. Raigorodskii, \emph{An improved upper bound on the chromatic number of the $n$-dimensional sphere}, \href{https://doi.org/10.48550/arXiv.2610.32408}{arXiv}.

\bibitem{Ro63} C. Rogers, \emph{Covering a sphere with spheres}, Mathematika 10 (1963), 157--164.
\bibitem{Ver} J.-L. Verger-Gaugry, {\it Covering a Ball with Smaller Equal Balls in $ {\mathbb R}^n $}, Discrete and Computational Geometry 33 (2005), no. 1, 143--155.
\bibitem{BW} K. B\"or\"oczky, G. Wintsche, {\it Covering the sphere by equal spherical balls}, Discrete and Computational
Geometry, The Goodman--Pollack Festschrift (2003), 237--253.

\bibitem{Jung} H. Jung, \emph{\"Uber die kleinste Kugel, die eine r\"aumliche Figur einschliesst}, Journal f\"ur die reine und angewandte Mathematik, 123 (1901), 241--257.

\bibitem{Pichugov} S. A. Pichugov, \emph{Jung's constant for the space $L_p$}, Mathematical notes of the Academy of Sciences of the USSR, 43 (1988), 604--614.

\bibitem{CGLW} F. Chen, S. Gao, X. Li, S. Wu, \emph{Covering the unit ball of $l^n_p$
with smaller balls and related
inequalities}, Math. Inequal. Appl. 28 (2025), no. 1, 99--119.

\bibitem{Rogers_density} C. A. Rogers, \emph{A note on coverings}, Mathematika 4 (1957), 1--6.

\bibitem{RZ} C. A. Rogers, C. Zong, \emph{Covering convex bodies by translates of convex bodies}, Mathematika, 44
(1997), no. 1, 215--218. 

\bibitem{RS} A. M. Raigorodskii, A. Sagdeev, \emph{A note on Borsuk's problem in Minkowski spaces}, Dokl. Math., 109 (2024), 80–83.

\bibitem{Melnikov} M. S. Mel'nikov, \emph{Dependence of volume and diameter of sets in an $n$-dimensional Banach space}, Uspekhi Matematicheskikh Nauk, 18 (1963), no. 4 (112), 165–170.

\bibitem{BeLi} Y. Benyamini, J. Lindenstrauss, \emph{Geometric Nonlinear Functional Analysis, Vol. 1}, AMS Colloquium Publications, 48 (2000).

\bibitem{Clarkson} J. A. Clarkson, \emph{Uniformly convex spaces}, Trans. Amer. Math. Soc., 40 (1936), 396–414.




\end{thebibliography}
\end{document}